\documentclass{amsart}   	
\usepackage{geometry}                		
\usepackage{graphicx,graphics, amsthm, amsfonts, amsbsy,amssymb, amsmath, cite, array, hyperref, enumerate}			
\renewenvironment{proof}{{\noindent\bfseries Proof.}}{\qed}								
\usepackage{amssymb}
\newtheorem{theorem}{Theorem}[section]

\newtheorem{corollary}[theorem]{Corollary}
\newtheorem{definition}[theorem]{Definition}
\newtheorem{proposition}[theorem]{Proposition}
\newtheorem{note}{Note}

\newtheorem{conjecture}{Conjecture}
\newtheorem{ass}{Assumption}

\newcommand{\grass}{\mathrm{Grass}}
\newcommand{\ext}{\mathrm{Ext}}
\newcommand{\expp}{\mathrm{Exp}}
\newcommand{\logg}{\mathrm{Log}}
\newcommand{\re}{\mathrm{Re}}
\newcommand{\im}{\mathrm{Im}}
\newcommand{\Arg}{\mathrm{Arg}}

\DeclareMathOperator{\tor}{Tor}

\author{Ralf Fr\"oberg}
\email{frobergralf@gmail.com}
\address{Department of Mathematics, Stockholm University, S-10691 Stockholm, Sweden}
\author{Clas L\"ofwall}
\email{clas.lofwall@gmail.com}
\address{Department of Mathematics, Stockholm University, S-10691 Stockholm, Sweden}

\usepackage{xcolor}

\title{Generic forms and algebras}

\begin{document}
\maketitle
\begin{abstract}
We begin with discussing what the concept of generic algebra should mean. In any case it should mean that it has minimal Hilbert series in
a sense we make precise.    In the commutative case there is a conjecture on what this minimal series is, and we give a conjecture for the generic series in the non-commutative quadratic case (building on work by Anick). 
 Furthermore we prove that if
$A=k\langle x_1,\ldots,x_n\rangle/(f_1,\ldots,f_r)$, $\deg(f_i)=2$ for all $i$, is a generic algebra, then $\{ x_if_j\}$ either is linearly independent or generate $A_3$. This complements a
similar theorem by Hochster-Laksov in the commutative case.  

We also consider universal enveloping algebras of Lie algebras with minimal quadratic presentations, and give a conjecture for the Hilbert series in the generic case,
and prove the conjecture when $r\le n^2/4$ or $r\ge (n^2-1)/3$, where $n$ is the number of generators and $r$ is the number of relations. 

Finally we show, a bit to our surprise, that the Koszul dual of a generic Koszul algebra is not
generic in general. But if the relations have algebraically independent coefficients over the prime field, we prove that the Koszul dual is generic. 
Hereby, we give a counterexample to a proposition by Polischchuk and Politselski, which states a criterion for a generic non-commutative quadratic presentation to be Koszul. 
We formulate and prove a correct version of the proposition.
\end{abstract}

\section{Introduction}
Let $F$ be either the free associative algebra in $n$ variables over a field or the polynomial ring $k[x_1,\ldots,x_n]$.
A presentation  $F/(f_1,\ldots,f_r)$ of a graded algebra is said to have type $(n;d_1,\ldots,d_r)$ if $f_i$ is a form of degree $d_i\ge2$,
$i=1,\ldots,r$. 
In the commutative case an algebra was called generic in \cite{Fr} if all coefficients of the $f_i$'s were
algebraically independent over the prime field. It was shown that generic algebras had minimal series in its type
and a formula for the Hilbert series was conjectured. This formula has been shown to be true in many special
cases, the most prominent for three variables (Anick) and for $r=n+1$ (Stanley). The
definition of generic above is very restrictive, and more general definitions, which also included the
non-commutative case were introduced in \cite{Fr-L} (for which this paper is a continuation) and in \cite{An2}. In most articles dealing with generic algebras, the concept "generic" is not defined. 
We have also seen the definition "a point in a countable intersection of non-empty Zariski open sets of an appropriate
affine space", without referring to Hilbert series, as in \cite[page 123]{P-P}, which we find strange. 

We begin this article with
a discussion of the concept generic algebra.
 We have tried to lay a self contained foundation for the concept "generic algebra" in the next
section. This means that some results, such as Proposition \ref{min} might be well known.

\section{Generic algebras}

By an algebra we mean a non-negatively graded algebra over a field $k$, $A=\oplus_{i\ge0} A_i$, where $A_0=k$ and $A$ is generated by the finite dimensional
vector space $A_1$. We study both the commutative and non-commutative case. Let $F$ be the free associative algebra $k\langle x_1,\ldots,x_n\rangle$
or the polynomial ring $k[x_1,\ldots,x_n]$, where $n=\dim (A_1)$. An algebra $A$ is the quotient of $F$ by a homogeneous ideal $I \subseteq \oplus_{i\ge2} F_i$. 
If $I$ is generated by $f_1,\ldots,f_r $ and $\deg(f_i) =d_i$, where $2\le d_1\le d_2\le\cdots\le d_r$, we say that $A$ has a {\em presentation of type $t=(n;d_1,\ldots,d_r)$} 
and $(f_1,\ldots,f_r)$ is said to be a {\em  sequence of type $t$}. If $I$ is minimally generated by $f_1,\ldots,f_r $ we say that $A$ has type $t$ as an algebra.
The type of an algebra is determined by $n$ and the graded vector space $\tor^A_2(k,k)$. In case $d_i=2$ for all $i=1,\ldots,r$, we call the algebra quadratic, 
and denote its type  by $(n,r)$. In the quadratic case, the presentation is minimal if and only if the $f_i$'s are linearly independent. 

\medskip
What should be meant by a generic algebra of type $t$? It is natural to require that, if $A$ is  generic, then ``almost all" algebras of type $t$ should have
the same Hilbert series as $A$, so if we pick an algebra ``at random", it will have this Hilbert series with high probability.
When the field is the real or complex numbers our definition of generic given below will satisfy this requirement. 
In the rest of the paper $k$ will be a field and $k_0$  a fixed but arbitrary prime field.
If $A$ is a graded algebra we denote the Hilbert series of $A$ by $A(z)$.

\begin{definition}
Let $A(z)=\sum_{n\ge0}a_nz^n$ and $B(z)=\sum_{n\ge0}b_nz^n$ be power series with real coefficients. We write $A(z)\ge B(z)$ if $a_n\ge b_n$ for all $n\ge0$.
\end{definition}

We start by giving, for each $t$, a concrete algebra with minimal Hilbert series.
\begin{definition}
Let $t=(n;d_1,\ldots,d_r)$ be a type. 
Let  $M_i^n$ be the set of monomials of degree $i$ in $k\langle x_1,\ldots,x_n\rangle$, 
 and $CM_i^n$ be the set of monomials of degree $i$ in  $k[x_1,\ldots,x_n]$. 
 Let $c=\{c_m^i;\ m\in CM_{d_i}^n, i=1,\ldots,r\}$ be a set of indeterminates.
 Put 
$$
A_{t}=k_0(c)[x_1,\ldots,x_n]/(f_{1,c},\ldots,f_{r,c}) \text{ where } f_{i,c}=\sum_{m\in CM_{d_i}^n}c_{m}^im.
$$ 
Analogously, in the non-commutative case, let  $c=\{c_m^i;\ m\in M_{d_i}^n, i=1,\ldots,r\}$ be a set of indeterminates. Put 
$$
A_{t}=k_0(c)\langle x_1,\ldots,x_n\rangle/(f_{1,c},\ldots,f_{r,c})\text{ where }f_{i,c}=\sum_{m\in M_{d_i}^n}c_{m}^im.
$$

\end{definition}
\begin{proposition}\label{min}
The Hilbert series $A_{t}(z)$ is minimal among all series $A(z)$ where $A$ is a commutative (non-commutative respectively) algebra with a presentation of type t over a field $k$ with prime field $k_0$.  
\end{proposition}

\begin{proof} We prove the commutative case, the non-commutative case is similar. Let $k$ be any field with $k_0$ as prime field and let $A$
be any commutative algebra with a presentation (not necessarily minimal) $A=k[x_1,\ldots,x_n]/(f_1,\ldots,f_r)$ of type $t$.
Consider the algebra $A'$ defined by the presentation $A'=k(c)[x_1,\ldots,x_n]/(f_{1,c},\ldots,f_{r,c})$ of type $t$. Then $A$ is a specialization of $A'$. The dimension 
of $(f_1,\ldots,f_r)_d$
is the rank of the matrix with rows $\{m_jf_i\}$, where $m_j$ is a monomial of degree $d-\deg(f_i)$ and analogously for $(f_{1,c},\ldots,f_{r,c})_d$. 
But if a subdeterminant of the latter matrix is zero, then also the corresponding subdeterminant of the former is zero. 
 Thus $\dim(f_1,\ldots,f_r)_d\le\dim(f_{1,c},\ldots,f_{r,c})_d$, and hence $\dim A_d\ge\dim A'_d$, so $A(z)\ge A'(z)$. Moreover $A'=A_t\otimes_{k_0}k$, and hence $A'(z)=A_t(z)$.
\end{proof}

\begin{definition}\label{gdef}  A presentation of a commutative (or non-commutative) algebra $A$ of type $t$ is called generic if $A(z)=A_t(z)$.  An algebra is called generic if its minimal presentation is generic. An algebra is called strictly generic if its minimal presentation has the property that all coefficients in the relations are algebraically independent over the prime field. 
 If $A$ is generic with a presentation $F/(f_1,\ldots,f_r)$, then the sequence  $(f_1,\ldots,f_r)$ is said to be generic. 
\end{definition}

It is not true in general  that $A_t$ has type $t$ as an algebra. An obvious such case is $t=(2;2,2,2,2)$, since there are just 3 linearly independent forms of degree 2 in 2 variables. A more subtle case is $t=(2;2,2,3)$. In this case $A_t$ as an algebra has type $(2;2,2)$, since 2 general forms in degree 2 generate everything in degree 3. 
The following proposition explains the situation.
\begin{proposition}\label{notgen}
Suppose $t=(n;d_1,\ldots,d_r)$ is a type such that $A_t$ has not type $t$ as an algebra. Then there is $s<r$ such that $A_t$ has type $(n;d_1,\ldots,d_s)$ and $(A_t)_{d_{s+1}}=0$.
\end{proposition}
\begin{proof}
Suppose $t=(n;d_1,\ldots,d_r)$ and the type of $A_t$ is not $t$. Recall that we always assume $2\le d_1\le d_2\le\cdots\le d_r$. Let $(f_1,\ldots,f_r)$ be the relations in $A_t$. Since the relations are not minimal there is a minimal $s$ such that $f_{s+1}$ belongs to the ideal $(f_1,\ldots,f_s)$. Since $f_{s+1}$ is general, it follows that $(A_{t})_{d_{s+1}}=0$ and $A_t$ has type $(n;d_1,\ldots,d_s)$.
\end{proof}

\medskip
We will often assume that $A_t$ has type $t$ as an algebra. This condition is true if $A$ is strictly generic of type $t$. It is also true if $t=(n;d,d,\ldots,d)$ with $r$ copies of $d\ge2$ where $r\le{n+d-1\choose d}$. In this case the $r$ relations in $A_t$ of degree $d$ are linearly independent and hence  $A_t$ has type $t$ as an algebra.
In particular, the condition is true in the quadratic case.

\begin{proposition}\label{trans}
Assume  $t$ is a type and $A_t$ as an algebra has type $t$. Then $A$ is generic of type $t$ if and only if $A$ has type $t$ and $A(z)$ is minimal among all series $B(z)$, where $B$ is a commutative (non-commutative, respectively) $k$-algebra of type $t$ and $k$ has prime field $k_0$. 
If the field $k$ has sufficiently high degree of transcendency over the prime field  $k_0$ then there exists a commutative 
(non-commutative, respectively) algebra over $k$  which is generic of type $t$.
\end{proposition}
\begin{proof} We prove the commutative case. Suppose $A$ is generic of type $t$, i.e., $A(z)=A_t(z)$. If $B$ is an algebra of type $t$, then in particular $B$ has a presentation of type $t$ and the result follows from Proposition \ref{min}. Suppose on the other hand  that $A$ has type $t$ and $A(z)$  is minimal among all series $B(z)$, where $B$ is a commutative $k$-algebra of type $t$ and $k$ has prime field $k_0$. Since $A_t$ is an algebra of type $t$, it follows that $A(z)\le A_t(z)$ but by Proposition \ref{min} we also have $A_t(z)\le A(z)$, hence $A$ is generic of type $t$. 

Suppose now $g=\{g_m^i;\ m\in CM_{d_i}^n, i=1,\ldots,r\}$ is a set of algebraically independent elements over $k_0$ in $k$. 
We have a map of fields $k_0(c)\to k$, which maps $c$ to $g$. Let $A=A_t\otimes_{k_0(c)}k$. Then $A(z)=A_t(z)$ and  $A$ has type $t$ since by assumption $A_t$ has type $t$. Hence $A$ is generic of type $t$.
\end{proof}

\medskip
\begin{ass}

In the rest of the paper we will assume, unless otherwise stated, that the field $k$ has sufficiently high transcendence degree over its prime field  $k_0$, 
in particular $k$ may be the real or complex numbers, if $k_0$ is the rational numbers.
\end{ass}

\medskip

\medskip
We will also study universal enveloping algebras of Lie algebras, more precisely we will study Lie superalgebras (with squares of odd elements) generated by odd generators 
of degree one with relations of degree two. The universal enveloping algebras of such Lie algebras are exactly the Koszul duals of commutative quadratic algebras. In the sequel 
we will interpret  $[a,b]$  in an associative algebra as $ab+ba$ if $a,b$ are of odd degree.

\begin{definition} A $k$-Lie-algebra is said to have type $(n,r)$ if it is the free Lie $k$-superalgebra (with squares of odd elements) $\mathcal L_k(x_1,\ldots, x_n)$ on the 
odd elements $x_1,\ldots,x_n$  modulo the Lie ideal minimally generated by linearly independent elements $f_1,\ldots,f_r$, where 
$f_i=\sum_{1\le j\le n}c_{i,j}x_j^2+\sum_{1\le j<l\le n}c_{i,j,l}[x_j,x_l]$ and $c_{i,j}, c_{i,j,l}\in k$. If $f_1,\ldots,f_r$ are interpreted as elements in  $k\langle x_1,\ldots,x_n\rangle$, 
then $k\langle x_1,\ldots,x_n\rangle/(f_1,\ldots,f_r)$ is said to have Lie type $(n,r)$.

Let $k_0$ be the prime field of $k$ and let $t=(n,r)$, $r\le {n+1\choose2}$, 
be a given type. Let $c=\{c_{i,j},1\le i\le r,1\le j\le n\}\cup\{c_{i,j,l},1\le i\le r,1\le j<l\le n\}$ be a set of indeterminates.
 Put 
$$
L_t=\mathcal L_{k_0(c)}(x_1,\ldots, x_n)/I_c\text{ and } A_t=k_0(c)\langle x_1,\ldots,x_n\rangle/I_c, 
$$ 
where $I_c=(f_{1,c},\ldots,f_{r,c})$ and 
$f_{i,c}=\sum_{1\le j\le n} c_{i,j}x_j^2+\sum_{1\le j<l\le n}c_{i,j,l}[x_j,x_l]$. 
\end{definition}

\begin{proposition}\label{lie} The Hilbert series  $L_t(z)$ is minimal among all series $L(z)$ where $L$ is a  $k$-Lie-algebra of type $(n,r)$ 
and $k$ has prime field $k_0$ and the Hilbert series $A_t(z)$ is minimal among all series $C(z)$ where $C$ is 
a $k$-algebra of Lie type $(n,r)$ and $k$ has prime field $k_0$. Also $UL_t=A_t$. Moreover, if $C=UL$ is any $k$-algebra of Lie type $(n,r)$ and $k$ has prime field $k_0$, then $C(z)=A_t(z)$ if and only if $L(z)=L_t(z)$. 
Also,  there exists a $k$-algebra $C$ of Lie type $(n,r)$ such that $C(z)=A_t(z)$ and there exists a $k$-Lie-algebra $L$ of type $(n,r)$ such that $L(z)=L_t(z)$.  
\end{proposition}

\begin{proof} Similar to the proofs of Proposition \ref{min} and Proposition \ref{trans}. Observe that the condition $A_t$ has type $t$ is satisfied since the relations are quadratic. The Hilbert series $UL(z)$ 
 is obtained from $L(z)=\sum_{i=1}^\infty e_iz^i$ by the exponential operator
\begin{align}
\expp(L)(z)=\prod_{i=1}^\infty\frac{(1+z^{2i-1})^{e_{2i-1}}}{(1-z^{2i})^{e_{2i}}}.\label{exp}
\end{align}
 This operator  has an inverse 
\begin{align}
\logg(V)(z)=\sum_{r=1}^\infty\frac{\mu(r)}{r}\log(V((-1)^{r+1}z^r))\label{log}
\end{align}
where $V(z)$ is a series with constant term 1, $\mu$ is the Möbius function and log is the natural logarithm, see \cite{Lof2}. If $L_1(z)=L_2(z)$ then the product 
formula (\ref{exp}) gives $UL_1(z)=UL_2(z)$. Also $L(z)$ is determined from $UL(z)$ by (\ref{log}). Hence the claim in the theorem that 
$C(z)=A_t(z)\iff L(z)=L_t(z)$ follows.
\end{proof}

\begin{definition}
A $k$-algebra $C$ of a certain Lie type  $t=(n,r)$ (a $k$-Lie-algebra $L$ of a certain type $t=(n,r)$, respectively), where  $r\le {n+1\choose2}$, is called generic if $C(z)=A_t(z)$ ($L(z)=L_t(z)$).
\end{definition}
\begin{proposition}\label{qgen}
Suppose $A$ is a $k$-algebra (commutative, non-commutative or of Lie type) of type $t=(n,r)$. Then $A$ is generic if and only if 
$A(z)$ is minimal among all series for $k$-algebras of type $t$.
\end{proposition}
\begin{proof} We know that $A_t$ has type $t$ as an algebra since $A$ is quadratic. Hence the result follows from Proposition \ref{trans}, which is true also in the Lie case.
\end{proof}

\vspace{12pt}
Let $t=(n;d_1,\ldots,d_r)$ be a type. Given a total order of the monomials, a point in $k^{n^d}$ corresponds to a polynomial of degree $d$ in $k\langle x_1,\ldots,x_n\rangle$.  Thus the points in $k^N$, 
$N=\sum_{i=1}^rn^{d_i}$ correspond to non-commutative presentations of type $t$.
Analogously the points in $k^N$, $N=\sum_{i=1}^r{n+d_i-1\choose n-1}$ correspond to commutative presentations of type $t$,
and in the Lie case, the points in $k^N$, $N=r{n+1\choose 2}$, correspond to Lie-algebras with $r$ quadratic relations. 
 
\begin{theorem}\label{open}
In any of the three cases above, consider presentations of a type $t$. Call the points which give generic algebras of type $t$
generic points. The set of generic points of the three classes in 
the appropriate space $k^N$  is an intersection of countably many Zariski-open sets.  The set is Zariski-open in the commutative case. It is also 
Zariski-open for non-commutative and Lie algebras, provided there are artinian examples (nilpotent, respectively) of the corresponding type.

The set of generic points is non-empty for the Lie case and it is also non-empty for the first two classes if $A_t$ as an algebra has type $t$.

Moreover, if $k$ is equal to the real or complex numbers then the set of points with Hilbert series not equal to the minimal series is a set of Lebesgue measure zero. 
\end{theorem}

\medskip
\begin{proof}
Let $\sum g_jz^j$ be the minimal series (defined in Proposition \ref{min} and \ref{lie}) in any of the three classes of algebras of type $t$.  Let $I$ be the ideal defined by a point. Let $F_{i,m}\subset k^N$
be the subspace of points $P$ with $\dim I_i<m$. Since $\dim I_i$ equals the rank of a large matrix, it follows that $F_{i,m}$ is closed, being the zero set of determinants.
Hence the set of points $F_i$ with Hilbert series $\sum \lambda_j(P)z^j$ and $\lambda_i(P)>g_i$ is closed,
and $\cup_{i=0}^\infty F_i$ consists of all points $P$ such that the corresponding algebra has not the minimal Hilbert series. The complement of $\cup_{i=0}^\infty F_i$ is an intersection of 
Zariski-open sets. The set of points which define a minimal presentation is Zariski-open. Hence the set of generic  points is an intersection of 
Zariski-open sets and it is non-empty, if $A_t$ has type $t$ as an algebra, by Proposition \ref{trans} and Proposition \ref{lie}. For the commutative case, it is shown in \cite[Proposition 1]{F-G-L}
that there are only finitely many possibilities for the Hilbert series of algebras of a given type. Since the Hilbert function eventually is a polynomial, there is a $C$, such that one of these
series has smallest coefficient of $z^i$ for all $i\ge C$. This gives that $\cup_{i=0}^\infty F_i=\cup_{i=0}^C F_i$.
We get the conclusion, that the set of algebras with minimal series  is open. Now consider the non-commutative or Lie case and
suppose there are artinian (nilpotent) examples. For a proof in this case, see also \cite[Appendix A, Theorem A.1]{L-N}. Then there exists a $C$ such that $g_C=0$. 
Hence if $P\in F_i$ for some $i>C$ then 
$\lambda_i(P)>0$ and it follows that also $\lambda_C(P)>0$, i.e., $P\in F_C$. Hence $\cup_{i=0}^\infty F_i=\cup_{i=0}^C F_i$ and it follows that the complement is open. 

In all cases, 
if $k$ is equal to the real or complex numbers it is shown in \cite{Ca} that an algebraic set has Lebesgue measure 0, so each $F_i$ has Lebesgue measure 0, 
and hence also $\cup_{i=0}^\infty F_i$ has measure 0.
\end{proof}

\medskip
\begin{theorem}\label{gendef} Suppose $t$ is a type such that $A_t$ has type $t$.
A commutative algebra of type $t$ is generic if and only if it belongs to a non-empty Zariski-open subset of the affine space of commutative presentations of type $t$, whose  
points have constant Hilbert series. A non-commutative or Lie algebra over the real or complex numbers of type $t$ is generic if and only if it belongs to a non-empty 
countable intersection of Zariski-open subset of the affine space of non-commutative or Lie algebras of type $t$, whose  points have constant Hilbert series.
\end{theorem}

\begin{proof}
If the algebra is generic of type $t$, the result follows from Theorem \ref{open}. Suppose $A$ is a commutative algebra of type $t$, which belongs to a non-empty Zariski-open subset $O$ of 
the affine space of commutative algebra presentations of type $t$ with the property that all points in $O$ have the same Hilbert series. By Theorem \ref{open}, the set $G$ of all generic points 
is open and non-empty. Since the field $k$ is infinite, $O\cap G$ is non-empty and from this it follows that $A$ is generic. Suppose now that $A$ is a non-commutative or Lie algebra 
of type $t$ over the real or complex numbers, which belongs to a non-empty countable intersection $O$ of Zariski-open subsets of the affine space of non-commutative or 
Lie algebra presentations of type $t$ with the property that all points in $O$ have the same Hilbert series. By Theorem \ref{open}, the set $G$ of all generic points is a non-empty countable 
intersection of open sets. By \cite{Ca} it follows that $G\cap O$ is non-empty since its complement is a countable union of sets of Lebesgue measure zero. 
We may conclude as before that $A$ is generic.
\end{proof}

\section{Minimal series}
\begin{definition}
For a power series $A(z)=\sum_{i\ge0}a_iz^i$, let $[A(z)]=\sum_{i\ge0}b_iz^i$, where $b_i=a_i$ if $a_j>0$ for all $j\le i$ and $b_i=0$ otherwise.\\
Also $\sum_{i\ge0}a_iz^i>_{lex}\sum_{i\ge0}b_iz^i$ if for some $n$ we have $a_i=b_i$ if $i<n$ and $a_n>b_n$.
For any type $t=(n;d_1,\ldots,d_r)$ let $F_t(z)=\left[\prod_{i=1}^r(1-z^{d_i})/(1-z)^n\right]$.

\end{definition}

\begin{conjecture}\label{conj}(\cite{Fr}) If $A$ is a commutative strictly generic algebra of type $t$ then $A$ has Hilbert series 
$F_t(z)$.
\end{conjecture}
It is shown in \cite{Fr} that $A_t(z)\ge_{lex} F_t(z)$. Also see \cite[Proposition 1]{Di}
\begin{note}
Let $I$ be a graded ideal in $S=k[x_1,\ldots,x_n]$ and let $A=S/I$. In \cite{Pa} a form $f\in S_d$ is called semiregular on $A$ if the multiplication $f\colon A_{e-d}\rightarrow A_e$
has full rank (is injective or surjective) for all $e$. A sequence $f_1,\ldots,f_r$ in $S$, $\deg(f_i)=d_i$, is called semiregular if the multiplication 
$f_i\colon S/(f_1,\ldots,f_{i-1})\rightarrow S/(f_1,\ldots,f_{i-1})$ is semiregular for all $i$. It is shown that $S/(f_1,\ldots.f_r)(z)=F_t(z)$ if and only if
$(f_1,\ldots,f_r)$ is semiregular.
 
 In \cite{Di} an element $f\in S_d$ is said to be regular up to gegree $D$ on $S/I$ if the multiplication $f\colon (S/I)_{e-d}\rightarrow (S/I)_e$ is injective for all $e\le D$.
 It is shown that if $S/(f_1,\ldots,f_r)$ is artinian, then $S/(f_1,\ldots.f_r)(z)=F_t(z)$ if and only if the multiplication 
 $f_i\colon S/(f_1,\ldots,f_{i-1})\rightarrow S/(f_1,\ldots,f_{i-1})$ has bounded regularity $D$ for all $i$, where $D$ is the highest nonzero degree in $S/(f_1,\ldots,f_r)$.
\end{note}

Conjecture \ref{conj} is proved for $r\le n$, \cite{Sta}, $n=2$, \cite{Fr}, $n=3$, \cite{An1}, $r=n+1$, \cite{St}, and for many other special commutative types \cite{Fr-Ho,Ne,Ni}.
In every case, the method of proof has been the following. 
If $A$ and $A_t$ has type $t$, we have $A(z)\ge A_t(z)\ge_{lex} F_t(z)$. Thus, if we find an algebra $A$ of type $t$ with $A(z)=F_t(z)$, we have proved Conjecture 1 for this type. This is the
method used in the cases when the conjecture is proved.

\begin{proposition} If Conjecture \ref{conj} is true, then the formula $A(z)=F_t(z)$ holds for all generic algebras $A$ of type $t$.
\end{proposition}
\begin{proof}
 By  assumption, $A_t(z)=F_t(z)$ if $A_t$ has type $t$. Let $t=(n;d_1,\ldots,d_r)$. By Proposition \ref{notgen} there is $s<r$ such that $A_t$ has  type $t'=(n;d_1,\ldots,d_s)$ and $(A_t)_{d_{s+1}}=0$. Since $A_{t'}$ has type $t'$, it follows by assumption that $A_{t'}(z)=F_{t'}(z)$, which is a polynomial of degree $<d_{s+1}$. The factors $(1-z^{d_i})$ for $i>s$ has no effect on $F_{t'}(z)$. Thus $F_t(z)=F_{t'}(z)$ and  hence   $A_t(z)=A_{t'}(z)=F_t(z)$ from which the result follows. 
\end{proof}

\vspace{15 pt}
For the non-commutative case, we need the notion of ``strongly free" defined by Anick, see \cite{An} and \cite{An2}. We choose the following definition. 
\begin{definition} A sequence $(f_1,\ldots,f_r)$ in $k\langle x_1,\ldots,x_n\rangle$ is called strongly free and also the quotient algebra $B$ is called strongly free 
if $B$ has global dimension 2.
\end{definition}

Let $t=(n;d_1,\ldots,d_r)$ 
be a type and let $B=k\langle x_1,\ldots,x_n\rangle/(f_1,\ldots,f_r)$ be an algebra of type $t$. There is an exact sequence of left $B$-modules, which is the beginning of a minimal resolution of $k$ as  a left module over $B$, 
\begin{align}\label{free}
\bigoplus_{i=1}^rBS_i\xrightarrow{\partial_2}\bigoplus_{j=1}^nBT_j\xrightarrow{\partial_1}B\to k\to0
\end{align}
where $S_i$ is a basis element of degree $d_i$ for $i=1,\ldots,r$ and $T_j$ is a basis element of degree 1 for $j=1,\ldots,n$. Write $f_i=\sum_{j=1}^nf_{ij}x_j$, where $f_{ij}$ 
are polynomials of degree $d_i-1$ for $i=1,\ldots,r$, $j=1,\ldots,n$. The map $\partial_2$ is given by $\partial_2(S_i)=\sum_{j=1}^nf_{ij}T_j$ for $i=1,\ldots,r$ and the map $\partial_1$ is given by $\partial_1(T_j)=x_j$ for $j=1,\ldots,n$.

It follows that an algebra $B$ is strongly free if and only if the map $\partial_2$ in (\ref{free}) is injective. It also follows that a reordering of  a strongly free sequence is still strongly free and that a subsequence of a strongly free sequence is strongly free. 

\begin{proposition}\label{eqfree}
Let $p_t(z)=1-nz+\sum_{i=1}^rz^{d_i}$ for a given type $t=(n;d_1,\ldots,d_r)$ and let $B$ be an algebra of type $t$.
\begin{align*}
(i)&\ B(z)p_t(z)\ge1\text{ with equality  if and only if $B$ is strongly free. } \\
(ii)&\ B(z)\ge_{lex}p_t(z)^{-1}.\\
(iii)&\ B(z)\ge_{lex} [p_t(z)^{-1}].
\end{align*} 
\end{proposition}

\begin{proof}
$(i)$ Let $K$ be the kernel of $\partial_2$. Taking the alternating sum of dimensions in different degrees of the exact sequence (\ref{free}) we get 
$B(z)p_t(z)=1+K(z)\ge1$.

$(ii)$ We use the following fact. If $R(z)\ge_{lex}S(z)$ and $T(z)\ge_{lex}0$, then $R(z)T(z)\ge_{lex}S(z)T(z)$. Since the constant term of $p_t(z)^{-1}$ is 1 we have $p_t(z)^{-1}>_{lex}0$ and since by $(i)$, $B(z)p_t(z)\ge_{lex}1$, we get $B(z)\ge_{lex}p_t(z)^{-1}$.

$(iii)$ This follows from $(ii)$.
\end{proof}

\medskip
In fact, the stronger inequality $B(z)\ge p_t(z)^{-1}$ is proven in \cite[Lemma 6.1]{An2}. However, we haven't been able to show the corresponding strong inequality in the commutative case.

\medskip
Unfortunately, it is not true that the generic series of type  $t$ is $[p_t(z)^{-1}]$. A counterexample is given in \cite[Example 6.3]{An2} with $t=(2;2,7)$. However, in the quadratic case we make the following conjecture.

\begin{conjecture} \label{art}If $B$ is a non-commutative generic algebra of type $(n,r)$ then $B$ has Hilbert series $[1/(1-nz+rz^2)]$. 
\end{conjecture}

Compare \cite[Question 6.4]{An2}. 
If the conjecture is true, then it follows that $B$ is artinian if and only if $r>n^2/4$, since the series $1/(1-nz+rz^2)$ has negative terms if $r>n^2/4$ and only positive terms if $r\le n^2/4$. 
This last fact is well known (see \cite[Theorem 5.1, Lemma 5.10]{An2}), but as convenience for the reader we give a proof of it here.  

\begin{proposition}\label{neg}
The power series expansion of $1/(1-nz+rz^2)$ has all coefficients positive if $r\le n^2/4$ and has some negative coefficient if $r>n^2/4$.
\end{proposition}

\begin{proof}
We have $1-nz+rz^2=(1-\frac{n}{2}z)^2-(\frac{n^2}{4}-r)z^2$. Suppose $r\le n^2/4$. Put $b=\sqrt{\frac{n^2}{4}-r}$. Then $1-nz+rz^2=(1-(\frac{n}{2}+b)z)(1-(\frac{n}{2}-b)z)$. 
We have $b\le n/2$. If $b=n/2$, then $r=0$ and the series expansion of $1/(1-nz)$ has only positive coefficients. If $b<n/2$ the series expansions of $1/(1-(\frac{n}{2}\pm b)z)$ 
have only positive coefficients and hence the same is true for $1/(1-nz+rz^2)$.

Next suppose $r>n^2/4$ and put $b=\sqrt{r-\frac{n^2}{4}}$ and $a=\frac{n}{2}+ib$. Then  $1-nz+rz^2=(1-az)(1-\bar az)$.  We have 
$1/(1-nz+rz^2)=A/(1-az)+\bar A/(1-\bar az)$, where $A=a/(a-\bar a)$. Hence $1/(1-nz+rz^2)=\sum_{k\ge0}(Aa^k+\bar A{\bar a}^k)z^k$ and
$$
Aa^k+\bar A{\bar a}^k=2\re(Aa^k)=2\re(\frac{1}{a-\bar a}a^{k+1})=\frac{1}{b}\re(\frac{1}{i}a^{k+1})=\frac{1}{b}\im(a^{k+1}).
$$ 
Since $b>0$ we have $0<\Arg(a)<\pi/2$ and hence $\pi<\Arg(a^{k+1})<2\pi$ for some $k$. 
\end{proof}

\begin{corollary} Suppose Conjecture \ref{art} is true and $B$ is a non-commutative generic algebra of type $(n,r)$. Then $B$ is artinian iff $r>n^2/4$ and the following are equivalent for $k\ge2$
\begin{align*}
(i)&\ B_k=0\\
(ii)&\ k\ge\pi/\arctan\sqrt{4r/n^2-1}-1\\
(iii)&\ r\ge(\tan^2(\pi/(k+1))+1)n^2/4.
\end{align*}
In particular, $B_3=0$ iff $r\ge n^2/2$, $B_4=0$ iff $r\ge \frac{1}{2}(3-\sqrt5)n^2$ and $B_5=0$ iff $r\ge n^2/3$.
\end{corollary}
\begin{proof} The first claim follows from Proposition \ref{neg}. Clearly $(ii)$ and $(iii)$ are equivalent. Suppose $(ii)$ is true. We will use the proof of Proposition \ref{neg}.  Since 
$\Arg(a)=\arctan\sqrt{4r/n^2-1}$ we have $\Arg(a^{k+1})\ge\pi$. Let $k_0$ be the least $k$ such that this is true. Then $B_{k_0}=0$ and hence also $B_k=0$. Suppose $B_k=0$. Let $k_0$ be the least $k$ such that $B_{k_0}=0$. We get $\Arg(a^{k_0+1})\ge\pi$ and hence $(k_0+1)\arctan\sqrt{4r/n^2-1}\ge\pi$ or $k_0\ge\pi/\arctan\sqrt{4r/n^2-1}-1$. Since $k\ge k_0$ $(ii)$ follows.
The other claims follow from $(iii)$, since $\tan(\pi/4)=1$, $\tan(\pi/5)=\sqrt{5-2\sqrt5}$ and $\tan(\pi/6)=1/\sqrt3$.
\end{proof}

\begin{note} Anick proves that $B_3=0$  for a generic algebra of type $(n,r)$ where $r\ge n^2/2$. See \cite[Example 6.2]{An2} and also Theorem \ref{strong} and Theorem \ref{indep}.
\end{note}

\begin{theorem}\label{strong}  Conjecture \ref{art} is true for $r\le n^2/4$ and $r\ge n^2/2$ and for $n\le6$. 
\end{theorem}

\begin{proof}
To prove the conjecture, it suffices by Proposition \ref{eqfree} $(ii)$ to give examples of algebras of type $(n,r)$ with the conjectured Hilbert series. Suppose $r\le n^2/4$. By Proposition \ref{neg} and 
Proposition \ref{eqfree} it is enough to prove that there is a strongly free algebra of type $(n,r)$. As was noted above, a subset of a strongly 
free set is strongly free (see also  \cite[Lemma 2.7]{An}). Hence, it is enough to give an example of a strongly free algebra when $r=\lfloor n^2/4\rfloor$.
If $r=n^2/4$, $n$ even,
consider 
$$
B=k\langle x_1,\ldots,x_{n/2},y_1,\dots,y_{n/2}\rangle/(y_ix_j, 1\le i,j\le n/2).
$$ 
Then 
$B$ is isomorphic as a vector space to $k\langle x_1,\ldots,x_{n/2}\rangle\otimes k\langle y_1,\ldots,y_{n/2}\rangle$, which has the series $1/(1-nz/2)^2=1/(1-nz+rz^2)$ and 
it follows by Proposition \ref{eqfree} that $B$ is strongly free. 

If $n=2s+1$ is odd we let the variables be $x_1,\ldots,x_s,y_1,\ldots,y_{s+1}$ and define $B$ as before. Then the number of relations in $B$ is $r=s^2+s=\lfloor n^2/4\rfloor$. 
The series for $B$ is 
$$
B(z)=\frac{1}{1-sz}\cdot\frac{1}{1-sz-z}=1/(1-nz+rz^2)
$$ 
and hence $B$ is strongly free.

Now, suppose $r\ge n^2/2$.
Since $1/(1-nz+rz^2)=1+nz +(n^2-r)z^2+(n^3-2nr)z^3+\cdots$ it follows that 
$[1/(1-nz+rz^2)]=1+nz+(n^2-r)z^2$ if and only if $r\ge n^2/2$. To prove the conjecture for $r\ge n^2/2$, it is thus enough to give an example of an algebra $A$ of type 
$(n,r)$ such that $A_3=0$, since then a generic algebra $A$ of type $(n,r)$ also has $A_3=0$ and $A$ has linearly independent relations so its series is 
$1+nz+(n^2-r)z^2$. It is enough to find such an example for  $r=\lceil n^2/2\rceil$ and this is done in \cite[Example 6.2]{An2} (we will come back to this example in the proof of Theorem \ref{indep}).

A method to compute the Hilbert series of non-commutative 
graded algebras  is given in \cite{Lu}. 
Using this method, algebras of type $(n,r)$, where $n^2/4< r<n^2/2$, with the conjectured Hilbert series are given for $n\le6$ in \cite{Lu},
thereby proving the conjecture for $n\le6$.
\end{proof}

\begin{note} The best general result we know of concerning when non-commutative algebras of type $(n,r)$ are artinian is in \cite{Iy-Sh} where it is shown that 
$r>(n^2+n)/4$ generic forms implies that the algebra is artinian.\end{note}

\medskip
In the commutative case Conjecture \ref{conj} has been proved in the first non-trivial degree $\min\{ d_i\}+1$, \cite{Ho-La}, with generalizations in \cite{Au} and \cite{Mi-Mi}.
Thus, in particular, if $f_1,\ldots,f_r$ are generic quadratic forms in $A=k[x_1,\ldots,x_n]$, then it is proved in \cite{Ho-La} that $\{x_if_j\}$ are either linearly independent or generate 
$A_3$. We will prove the corresponding result in the non-commutative quadratic case. 

\begin{theorem}\label{indep}
Let $f_1,\ldots,f_r$ be generic quadratic forms in $A=k\langle x_1,\ldots,x_n\rangle$. Then $\{x_if_j\}\cup\{f_jx_i\},1\le i\le n,1\le j\le r$, are either linearly
independent or generate $A_3$.
\end{theorem}

\begin{proof}
For any quadratic forms $\{f_1,\ldots,f_r\}$ in $A$ we have that $\{x_if_j\}\cup\{f_jx_i\}$, $1\le i\le n,1\le j\le r$, generate a space of dimension $\le 2nr$ with equality if and only if
they are linearly independent.

Suppose $n=2s$ is even. Call the variables in $A$ $x_1,\ldots,x_s,y_1,\ldots,y_s$. There is in \cite[Example 6.2]{An2} an example with $n^2/2=2s^2$ forms $h_1(i,j)=x_iy_j$, 
$h_2(i,j)=x_ix_j-y_iy_j$, $1\le i,j\le s$, generating $A_3$ as a two-sided ideal. Since the number of elements of 
\begin{align}
\{x_ph_k(i,j),y_ph_k(i,j),h_k(i,j)x_p,h_k(i,j)y_p,1\le p,i,j\le s, k=1,2\}\label{even}
\end{align}
is $4\cdot 2s^3=n^3=\dim A_3$, we get that they are
linearly independent. Thus also,  if $f_1,\ldots,f_r$ are generic quadratic forms and $r\le n^2/2$, then $x_if_k,f_kx_i,y_if_k,f_ky_i$ are linearly independent. 
Since there are $n^2/2$ forms that generate $A_3$ as a two-sided ideal, also $r\ge n^2/2$
generic quadratic forms generate $A_3$ as a two-sided ideal.

Now suppose $n=2s+1$ is odd. Call the variables $x_1,\ldots,x_s,y_0,y_1,\ldots,y_s$. Then $n^2/2=2s^2+2s+1/2$. We will find $2s^2+2s$ quadratic forms $f_k$ in $2s+1$ variables such that 
$x_if_k,f_kx_i,y_if_k,f_ky_i$ are linearly independent (note that any subset of size $2s^2+2s$ of the $2s^2+2s+1$ forms given in \cite[Example 6.2]{An2} will  not satisfy this condition),
and $2s^2+2s+1$ quadratic forms $f_k$ in $2s+1$ variables such that $x_if_k,f_kx_i,y_if_k,f_ky_i$, generate $A_3$. 

We start with the $2s^2$ forms $h_1(i,j),h_2(i,j)$ from above of $y_0$-degree 0 and adjoin $2s$ new forms of $y_0$-degree 1. These are, for $1\le i\le s$,
\begin{align*}
q_1(i)&=y_0x_i+x_iy_0+y_iy_0\\
q_2(i)&=y_0y_i+y_iy_0.
\end{align*} 
Multiplying these $2s^2+2s$ forms to the right and to the left with all variables we get $(4s+2)(2s^2+2s)=8s^3+12s^2+4s$ expressions, which consist of three groups of 
$y_0$-degree 0,1,2, respectively.  The first group are the old $8s^3$
expressions (\ref{even}) of $y_0$-degree 0. The second group consist of $12s^2$ expressions which all have $y_0$-degree 1. These are, for $1\le i,j\le s$,
\begin{align*}
&y_0h_1(i,j),h_1(i,j)y_0,y_0h_2(i,j),h_2(i,j)y_0,\\
&x_iq_1(j),q_1(j)x_i,x_iq_2(j),q_2(j)x_i,y_iq_1(j),q_1(j)y_i,y_iq_2(j),q_2(j)y_i.
\end{align*}
Finally the third group are the following $4s$ expressions of $y_0$-degree 2,
$$y_0q_1(i),q_1(i)y_0,y_0q_2(i),q_2(i)y_0\text{ for }1\le i\le s.$$
To show that all these forms of degree 3 are linearly independent, we can treat the three groups separately. 

The first group we know are linearly independent from the even case above. 
Consider the second group. Suppose
\begin{align*}
\sum_{i,j=1}^s&\lambda_{1,i,j}y_0h_1(i,j)+\lambda_{2,i,j}h_1(i,j)y_0+\lambda_{3,i,j}y_0h_2(i,j)+\lambda_{4,i,j}h_2(i,j)y_0+\lambda_{5,i,j}
x_iq_1(j)+\\
&\lambda_{6,i,j}q_1(j)x_i+
\lambda_{7,i,j}x_iq_2(j)+\lambda_{8,i,j}q_2(j)x_i+\lambda_{9,i,j}y_iq_1(j)+\\
&\lambda_{10,i,j}q_1(j)y_i+\lambda_{11,i,j}y_iq_2(j)+\lambda_{12,i,j}q_2(j)y_i=0.
\end{align*}
The monomial $y_ix_jy_0$ occurs only in $y_iq_1(j)$, hence $\lambda_{9,i,j}=0$ for all $1\le i,j\le s$.
In the remaining expressions the monomial $y_0y_jx_i$ occurs only in $q_2(j)x_i$, hence $\lambda_{8,i,j}=0$ for all $1\le i,j\le s$.
Continuing like this, considering in turn, $x_jy_0x_i$, $y_0x_ix_j$, $x_iy_0x_j$, $x_ix_jy_0$, $y_iy_jy_0$, $y_0y_jy_i$, $y_jy_0y_i$,
$y_0x_iy_j$, $x_iy_0y_j$, $x_iy_jy_0$, we get successively 
$\lambda_{6,i,j}=\lambda_{3,i,j}=\lambda_{5,i,j}=\lambda_{4,i,j}=\lambda_{11,i,j}=\lambda_{12,i,j}=\lambda_{10,i,j}=\lambda_{1,i,j}=\lambda_{7,i,j}=\lambda_{2,i,j}=0$ for all $1\le i,j\le s$.

To prove that the  third group are linearly independent, suppose 
$$
\sum_{i=1}^s\lambda_{1,i}y_0q_1(i)+\lambda_{2,i}q_1(i)y_0+\lambda_{3,i}y_0q_2(i)+\lambda_{4,i}q_2(i)y_0=0.
$$
The monomial $y_0^2x_i$ only occurs once for each $i$, hence $\lambda_{1,i}=0$ for all $i$. Looking at the monomials $x_iy_0^2$ and $y_0^2y_i$ we get in the same way that 
$\lambda_{2,i}=\lambda_{3,i}=0$ for all $i$. Finally also $\lambda_{4,i}=0$ for all $i$.

We have thus found $2s^2+2s$ quadratic forms $f_k$ such that $x_if_k,f_kx_i,y_if_k,f_ky_i$ are linearly independent.
It follows that the same is true for $r$ generic quadratic forms if $r\le 2s^2+2s$. 

\medskip
Now add $y_0^2$ to the $2s^2+2s$ forms above. We will prove that they generate $A_3$ as an ideal. The first two groups are the same and the third group is enlarged by $x_iy_0^2,y_iy_0^2,y_0^2x_i,y_0^2y_i$ for $1\le i\le s$. There is also a new group of $y_0$-degree 3 consisting of just $y_0^3$. 

By the even case the first group generate everything in $A_3$ of $y_0$-degree 0. 

The second group consists of $12s^2$ linearly independent
elements of $y_0$-degree 1, and $12s^2$ is also the dimension  of $A_3$ in this $y_0$-degree and hence the second group generates everything in $A_3$ of $y_0$-degree 1. 

The third group contains $y_0q_2(i)=y_0^2y_i+y_0y_iy_0$ for all $i$ and hence $y_0y_iy_0$ is generated for all $i$. But $y_0q_1(i)=y_0^2x_i+y_0x_iy_0+y_0y_iy_0$ and hence also $y_0x_iy_0$ is generated for all $i$. Thus the third group generates  everything in $A_3$ of $y_0$-degree 2. 

Finally, the new group consisting of $y_0^3$ is a basis for $A_3$ in $y_0$-degree 3.

We have thus found $2s^2+2s+1$ quadratic forms $f_k$ such that $\{f_k\}$ generate $A_3$ as a two-sided ideal. It follows that the same is true for $r$ generic quadratic forms if $r\ge 2s^2+2s+1$. 

\end{proof}

\medskip
There is a natural analogue to Conjecture \ref{art} for Lie algebras (compare \cite[Conjecture 7.6]{Fr-L2} ).

\begin{conjecture}\label{liecon}
A generic Lie algebra of type $(n,r)$ has Hilbert series  $\left[\logg(1/(1-nz+rz^2))\right]$, where $\logg$ is defined in Section 2 equation (\ref{log}).
\end{conjecture}

If the conjecture is true  then $L$ is nilpotent if $r>n^2/4$, since then $1/(1-nz+rz^2)$ has negative terms (see Proposition \ref{neg} and \cite[Theorem 5.1]{An2}) and 
then also $\logg(1/(1-nz+rz^2))$ has negative terms, since if $V(z)\ge0$ then also $\expp(V)(z)\ge0$.

We have the following  analogue of Theorem \ref{strong}.

\begin{theorem}
Conjecture \ref{liecon} is true for $r\le n^2/4$ and for $r\ge(n^2-1)/3$. Moreover $L_3=0$ if and only if $r\ge(n^2-1)/3$ for a generic Lie algebra of type $(n,r)$. 
\end{theorem}
\begin{proof}
For the first case,  we proceed as in the proof of Theorem \ref{strong}. The only difference is that we change the relations $x_iy_j$ to $[x_i,y_j]$ and the isomorphism as vector spaces used in the proof is now true even as algebras.  Hence $UL(z)=1/(1-nz+rz^2)$ which gives $L(z)=\logg(1/(1-nz+rz^2))$.

We now prove the claim that a generic Lie algebra of type $(n,r)$ has $L_3=0$ if and only if $r\ge(n^2-1)/3$. 

Suppose $r\ge (n^2-1)/3$ and suppose $f_1,\ldots,f_s$ are generic quadratic forms in $k[x_1,\ldots,x_n]$, where $s={n+1\choose2}-r\le {n+1\choose2}-(n^2-1)/3={n+2\choose3}/n$. 
By \cite{Ho-La} $\{x_if_j\}$ are linearly independent and by \cite[Theorem 2.6]{Lof} it follows that $L_3=0$ where $UL$ is the Koszul dual of $k[x_1,\ldots,x_n]/(f_1,\ldots,f_s)$. 
Hence there is an example of a Lie algebra $L$ of type $(n,r)$ such that $L_3=0$ which implies that $L_3=0$ in the generic case. 

Suppose on the other hand that $L_3=0$ for a Lie algebra of type $(n,r)$. Then by \cite[Theorem 2.6]{Lof} it follows that the dual relations $f_1,\ldots,f_s$ satisfy that  $\{x_if_j\}$ 
are linearly independent. But then $ns\le{n+2\choose3}$ from which it follows that $r={n+1\choose2}-s\ge(n^2-1)/3$.

Finally we prove that the conjecture is true when $r\ge (n^2-1)/3$. By the above we have $L(z)=nz+sz^2$ for a generic Lie algebra of type $(n,r)$ with 
$r\ge (n^2-1)/3$. We have to compute the exponent $e$ such that 
$$
1/(1-nz+rz^2)\equiv (1+z)^n(1+z^3)^e/(1-z^2)^s \text{ modulo } z^4
$$ 
A direct calculation gives $e=n(n^2-(3r+1))/3$ and thus 
$$
\logg(1/(1-nz+rz^2))=nz+sz^2+n(n^2-(3r+1))/3z^3+\cdots
$$ 
and in particular, 
$$
\left[\logg(1/(1-nz+rz^2))\right]=nz+sz^2\iff r\ge (n^2-1)/3
$$
and hence $L(z)=\left[\logg(1/(1-nz+rz^2))\right]$.

\end{proof}

\section{The Koszul dual}\label{dual}
We may consider a quadratic non-commutative algebra $A$ with $n$ generators and $r$ linearly independent quadratic relations as  a point in the Grassmannian $\grass(n^2,r)$. 
This is in particular useful when the Koszul dual is studied. If $U$ is a subspace of $V\otimes V$ then the Koszul dual of $A=T(V)/\langle U\rangle$ is defined as 
$A^!=T(V^*)/\langle U^0\rangle$ where $U^0=\{f\in V^*\otimes V^*;\  f(U)=0\}$.  It is known that $A^!$ is the subalgebra of $\ext_A(k,k)$ generated by 
$\ext_A^{1}(k,k)$, see \cite{Lof}. The correspondence between $A$ and  $A^!$ may be seen as a map $\grass(n^2,r)\to\grass(n^2,n^2-r)$, given by $U\mapsto U^0$. 
Moreover, by performing Gaussian elimination, a basis for $U^0$ may be obtained using rational expressions in terms of the coefficients of a basis for $U$. 
Thus the map  $A\to A^!$ is a map of projective manifolds and it is an isomorphism, since $(A^{!})^!=A$. 

If $A$ is commutative we may consider  $A$ as an element in $\grass({n+1\choose2},r)$. In this case $A^!$ as a quotient of $k\langle T_1,\ldots,T_n\rangle$ has relations which 
are linear combinations of $T_i^2$ and 
$[T_i,T_j]=T_iT_j+T_jT_i$. Hence $A^!$ is the enveloping algebra of a graded Lie superalgebra generated by $n$ odd elements of degree 1.  
The number of relations in $A^!$ is $n^2-{n\choose2}-r={n+1\choose2}-r$. Hence, when $A$ is commutative, the correspondence $A\leftrightarrow A^!$ may be seen as a map 
$\grass({n+1\choose2},r)\to\grass({n+1\choose2},{n+1\choose2}-r)$.

It is reasonable to believe that if a quadratic $k$-algebra is strictly generic (i.e., has algebraically independent coefficients in the relations), then the same is true for the Koszul dual, 
but we have no proof. However, we can prove the following fact. 
\begin{theorem}\label{kos} Let $B$ be a strictly generic quadratic $k$-algebra (commutative, non-commutative or of Lie type). Then $B^!$ is a generic $k$-algebra 
(of Lie type, non-commutative, or commutative, respectively). 
\end{theorem}
\begin{proof} We prove the non-commutative case, the other cases are similar. Let $t=(n,r)$ be a given type and let $A_t$ be the algebra introduced in Proposition \ref{min}. 
As a first step, we prove that $A_t^!$ is generic. Suppose $C$ is any $k$-algebra of the type of $A_t^!$ with prime field $k_0$. Then $C^!$ is a $k$-algebra of the type $t$ 
and hence $C^!$ is a specialization of $A_t$. The same specialization may be used to compute the Koszul dual of $C^!$, which is equal to $C$. 
Indeed,  $(C^!)^!$ is obtained by first computing $A_t^!$ and then specialize. It follows that $C(z)\ge A_t^!(z)$ and hence $A_t^!$ is generic by Proposition \ref{qgen}. 

Now let $B$ be any strictly generic $k$-algebra, where $k$ has prime field $k_0$, of type $t$. Then (with notations from Proposition \ref{min}) there is  a map of fields 
$k_0(c)\to k$ such that $B=A_t\otimes_{k_0(c)}k$. We have $B^!=A_t^!\otimes_{k_0(c)}k$ and hence $B^!(z)=A_t^!(z)$ and it follows that $B^!$ is generic.
\end{proof}

\medskip
The following theorem shows that it is essential in Theorem \ref{kos}  to assume that $A$ is strictly generic.

\begin{theorem}\label{ex} Suppose a quadratic $k$-algebra is either
\begin{align*}
(a)& \text{ a generic commutative algebra or}\\
(b)& \text{ a generic non-commutative algebra or}\\
(c)& \text{ a generic algebra of Lie type} 
 \end{align*}
 Then the  Koszul dual is in general not generic in $(a)$ as an algebra of Lie type, in $(b)$ as a non-commutative algebra and in $(c)$ as a commutative algebra.
 \end{theorem}
 \begin{proof} $(a)$ 
 
 \medskip
 Consider the algebra of commutative type $(4,6)$, 
 $$
 A=k[ a,b,c,d]/(a^2,b^2,c^2,ad-d^2,ac-bd,cd).
 $$ 
 It is easy to see that $A_3=0$, hence $A$ is generic. We have 
 $$
 A^!=k\langle a,b,c,d\rangle/([a,b],[b,c],[a,c]+[b,d],[a,d]+d^2). 
 $$
 If $A^!$ is generic, then it is strongly free with series $1/(1-4z+4z^2)$. Using e.g. the Macaulay2 \cite{M2} package ``GradedLieAlgebras" \cite{L-L}, 
 one can see that $L(z)=4z+6z^2+4z^3+7z^4+\cdots$ where $UL=A^!$. Developing $1/(1-4z+4z^2)$ as an infinite product, e.g. by using (\ref{log}) Section 2, we get 
 $$
 1/(1-4z+4z^2)=\frac{(1+z)^4}{(1-z^2)^6}\cdot\frac{(1+z^3)^4}{(1-z^4)^6}\cdots.
 $$
  Hence $L(z)$ is not minimal and hence the $k$-Lie-algebra $L$ is not generic of type $(4,4)$ and then, 
 by Proposition \ref{lie} also the associative algebra $A^!$ is not generic of Lie type $(4,4)$.
 
 \medskip
 \noindent $(b)$
 
 \medskip
 The algebra $A$ defined above may also be considered as a non-commutative algebra of type $(4,12)$, namely
 $$
 A=k\langle a,b,c,d\rangle/(a^2,b^2,c^2,ad-d^2,ac-bd,cd,ab-ba,ac-ca,ad-da,bc-cb,bd-db,cd-dc).
 $$
We still have $A_3=0$ and $A^!$ is the same, so this gives the required example. We will give another example studied in \cite[Lemma 4]{Fr-L}
$$
B=k\langle a,b,c,d\rangle/(a^2,b^2,d^2,ad,ca,cd,da,db,dc,ba,ab+ac,bc+cb,bc+c^2).
$$
The Koszul dual of $B$ is
$$
B^!=k\langle a,b,c,d\rangle/(ab-ac,bc-cb-c^2,bd).
$$
Again, it is fairly easy to see that $B_3=0$ and hence $B$ is generic. If $B^!$ were generic then, by Proposition \ref{neg} and Theorem \ref{strong}, its series would be $1/(1-4z+3z^2)$. 
However, In \cite[Lemma 4]{Fr-L} it is proved that $B^!(z)=1/(1-4z+3z^2-z^4)$.

\medskip
\noindent $(c)$

\medskip
Consider the following algebra of Lie type $(3,3)$.
$$
C=k\langle a,b,c\rangle/(b^2,[a,b]-c^2,[a,c]).
$$
Also in this case we can use the Macaulay2 package ``GradedLieAlgebras" to get that $L_3=0$, where $C=UL$. Hence $L_3$ is generic and by Proposition \ref{lie} 
it follows that also $C$ is generic. The Koszul dual of $C$ is
$$
C^!=k[ a,b,c]/(a^2,bc,ab+c^2).
$$
Since $C^!$ has commutative type $(3,3)$ the generic case is a complete intersection. However, $ab$ and $ab+c^2$ is $\ne0$ in $k[a,b,c]/(a^2,bc)$ but $ab(ab+c^2)=0$ in $k[a,b,c]/(a^2,bc)$. 
This example, under the name 3b, was studied in \cite{Ba-Fr}.		
\end{proof}

\begin{note} In the proof of \cite[Proposition 4.2]{P-P} it is used that the Koszul dual of a generic quadratic algebra is generic. The definition of generic in \cite{P-P} is that the algebra 
should belong to a countable intersection $C$ of non-empty Zariski open subsets of the appropriate affine space. We think that one has to add the assumption that the Hilbert series 
is constant in $C$, so by Theorem \ref{gendef}  it coincides with our definition. But even with this assumption it is not true that the Koszul dual of a generic algebra is generic, by 
Theorem \ref{ex}.  Indeed, the example provided in the proof of b) in Theorem \ref{ex} gives a counterexample also to \cite[Proposition 4.2]{P-P}, since $A$ satisfies their condition 
$s\le m^2/4$, $A$ is not Koszul since $A^!$ is not Koszul and $A$ is generic since $A_3=0$. We give in Theorem \ref{genkos} below a correct version of \cite[Proposition 4.2]{P-P}.
\end{note}

\begin{theorem}\label{genkos}
A strictly generic non-commutative algebra of type $(n,r)$ is Koszul if and only if $r\le n^2/4$ or $r\ge 3n^2/4$.
\end{theorem} 
\begin{proof}
Suppose first $A$ is a generic non-commutative algebra of type $(n,r)$. If $r\le n^2/4$ then, by Theorem \ref{strong}, Proposition \ref{neg} and Proposition \ref{eqfree}, 
$A$ is strongly free and thus has global dimension $\le2$. Hence $A$ is Koszul, since $\tor^A(k,k)$ lies on the diagonal.  If $n^2/2\le r< 3n^2/4$, then by \cite[Example 6.2]{An2}, 
$A_3=0$ and if $A$ is Koszul then $A^!$ has Hilbert series $1/(1-nz+sz^2)$, where $s=\dim(A_2)=n^2-r>n^2/4$. However, by Proposition \ref{neg}, the series $1/(1-nz+sz^2)$ 
has some negative coefficient. Thus we have proved that if $A$ is generic of type $(n,r)$, then $A$ is Koszul if $r\le n^2/4$ and not Koszul if $n^2/2\le r<3n^2/4$. 

Suppose now that $A$ is strictly generic of type $(n,r)$ and $n^2/4< r\le n^2/2$. Then by Theorem \ref{kos} $A^!$ is generic  and its type is $(n,s)$ where $n^2/2\le s<3n^2/4$. 
By the above we conclude that $A^!$ is not Koszul and hence also $A$ is not Koszul. Finally, suppose $A$ is strictly generic of type $(n,r)$ where $r\ge 3n^2/4$. Then by 
Theorem \ref{kos} $A^!$ is generic of type $(n,s)$ where $s\le n^2/4$. Hence by the above $A^!$ is Koszul and hence also $A$ is Koszul.
\end{proof}

\medskip

\medskip\noindent
There is no conflict of interest. All data is created by the authors. No funding was received.

\end{document}